# Bio-Inspired Strategies for Optimizing Radiation Therapy under Uncertainties

Keshav Kumar K.[1], Dr.NVSL Narasimham[2]

Research Scholar, Department of Mathematics, Jawaharlal Nehru Technological University, Hyderabad 500 085, India[1]

Associate Professor, , Department of Humanities and Mathematics, G. Narayanamma Institute of Technology and Science (for Women), Hyderabad 500 104, India[2]

∗Corresponding Author's Email: keshav.gnits@gmail.com

**ABSTRACT**

Radiation therapy is a critical component of cancer treatment. However, the delivery of radiation poses inherent challenges, particularly in minimizing radiation exposure to healthy organs surrounding the tumor site. One significant contributing factor to this challenge is the patient's respiration, which introduces uncertainties in the precise targeting of radiation. Managing these uncertainties during radiotherapy is essential to ensure effective tumor treatment while minimizing the adverse effects on healthy tissues. This research addresses the crucial objective of achieving a balanced dose distribution during radiation therapy under conditions of respiration uncertainty. To tackle this issue, we begin by developing a motion uncertainty model employing probability density functions that characterize breathing motion patterns. This model forms the foundation for our efforts to optimize radiation dose delivery. Next, we employ three bio-inspired optimization techniques: Cuckoo search optimization (CSO), flower pollination algorithm (FPA), and bat search Optimization (BSO). Our research evaluates the dose distribution in Gy on both the tumor and healthy organs by applying these bio-inspired optimization methods to identify the most effective approach. This research ultimately aids in refining the strategies used in radiation therapy planning under the challenging conditions posed by respiration uncertainty. Through the application of bio-inspired optimization techniques and a comprehensive evaluation of dose distribution, we seek to improve the precision and safety of radiation therapy, thereby advancing cancer treatment outcomes.

## 1. INTRODUCTION

Radiotherapy is a medical procedure that uses ionizing radiation sources such as protons, electrons, and high-energy particles to slow the growth of malignant growths [1]. It is of the highest priority to successfully align ionizing radiation beams with the 3-D shape of the tumor while protecting adjacent healthy tissue [2, 3]. Yet, this task becomes progressively more complex when addressing tumors located in the thorax and abdominal areas due to the inherent motion of the tumor during the treatment process [4]. This motion is primarily induced by quasi-periodic breathing patterns and is particularly significant for thorax tumors like those in the lungs and breast [5]. The constant motion of tumors during radiotherapy presents a significant problem because the tumor's exact position is not consistently known. Among the various sources of uncertainties, this review primarily focuses on the intrafractional respiratory motion, which results from the involuntary physiological process of respiration [6]. The organs located within the thoracic and

Keshav Kumar K., Dr.NVSL Narasimham(2023)

upper abdominal regions, including the liver, lungs, prostate, pancreas, esophagus, breast, and kidneys, undergo motion as a result of breathing [7]. This movement brings about substantial uncertainties in various aspects, including imaging, treatment planning, and the administration of radiotherapy for thoracic and abdominal conditions. Incorporating margins is a common practice to address uncertainties in tumor localization during radiotherapy. Nevertheless, the use of these margins amplifies the potential for radiotherapy-associated toxicity, as it extends the reach of radiation to normal tissues within the Planned Target Volume (PTV). Most of these side effects are attributed to uncertainties in tumor localization caused by breathing-induced motion and setup errors. Radiation oncologists must carefully balance the clinical benefits of treatment with the risks to the patient's long-term quality of life when determining radiation dosage. This trade-off, due to uncertainties in tumor localization, can also hinder the effectiveness of radiotherapy by preventing the delivery of the necessary dose escalation for effective treatment.

In this research, we take lung cancer treatment. Lung cancer continues to hold the unfortunate distinction of being the foremost cause of cancer-related fatalities, not only in the United States but also globally [8]. Its annual death toll nearly equals the combined mortality rates of prostate, breast, and colon cancer. A 2020 report focusing on lung cancer emphasizes that it is the most commonly diagnosed cancer and the primary contributor to cancer-related deaths in Canada. Globally, cancer rates are projected to double by 2050, with lung cancer being the most prominent [9]. Radiotherapy is employed in the treatment of more than half of all cancer patients [10]. Specifically, we are directing our attention to external beam radiotherapy, a method that utilizes a linear accelerator affixed to a revolving gantry to administer high-energy photon beams to the patient. Photon beams deposit energy as they traverse tissue, affecting both tumor cells and the healthy tissue in their path. To minimize damage to healthy cells, radiation is delivered from various angles, allowing each beam to deliver a small dose to healthy tissue while concentrating a high dose in the overlapping region centred on the tumor.

Our research aims to comprehend the impact of motion uncertainty on lung radiotherapy quality and establish a framework that generates solutions resistant to this uncertainty. We propose a bio-inspired optimization approach specifically tailored to address motion uncertainty and demonstrate its effectiveness. Our goal is to strike a balance between protecting healthy tissue and effectively treating the tumor, taking into account the presence of uncertainty.

## 2. Literature Survey

The literature survey on optimizing radiotherapy under uncertainties presents a rich tapestry of research aimed at improving the precision and effectiveness of this crucial medical treatment. In the research [11], the primary objective is evident: to address the complexities introduced by intrafraction motion by employing feedback control of the radiation dose administered. This innovative technique combines pre-treatment 4-D computed tomography (4DCT) imaging with intrafraction respiratory-motion surrogates to estimate the total given dosage and the predicted motion trajectory throughout treatment in real-time. The optimization of intensity-modulated radiotherapy (IMRT) plans under free-breathing conditions is a significant advancement. Notably, this study demonstrates that the proposed stochastic control approach not only reduces irradiated tissue volume compared to traditional internal target volume (ITV) treatment but also significantly cuts down treatment time without compromising dosimetric quality. It represents a promising avenue to enhance the efficiency of radiotherapy, particularly in scenarios where respiratory gating may be impractical or less efficient. In the study [12], the focus transitions to the domain of 4D multi-image-based (4DMIB) optimization, a field with the potential to bolster the resilience of scanned particle therapy in the presence of motion induced by respiration. The review underscores the pressing need for more comprehensive clinical evidence regarding the essentiality of 4DMIB optimization, particularly for conditions influenced by anatomical variations. Despite the wealth of research and technical insights in this domain, clinical investigations remain sparse, often constrained by methodological limitations such as limited patient cohorts and considerations related to motion dynamics. Nevertheless, the report acknowledges that robust 3D optimized plans appear to conform well to clinical tolerances, rendering them suitable for treating mobile targets using scanned particle therapy. The clinical urgency for the adoption of 4DMIB optimization, however, is noted to be contingent upon more substantial empirical demonstration.

In the study [13], the development of a risk-based robust approach is introduced, with a particular focus on addressing uncertainties related to tumor shrinkage during radiotherapy. The core objective of this suggested model is to reduce the variability of delivered doses, especially in worst-case scenarios, and minimize total radiation exposure to healthy tissues. The model leverages adaptive radiotherapy, a fractionation technique that considers the tumor's response to treatment over time and re-optimizes the treatment plan based on an estimate of tumor shrinkage. The clinical application of this approach is exemplified through a case study of lung cancer. The outcomes of this investigation highlight the potential benefits of the robust-adaptive model in terms of ensuring dose consistency within the tumor target while minimizing the impact on organs at risk. Furthermore, the model demonstrates superior performance in terms of maintaining uniform tumor dose distribution and overall plan reliability, underscoring its potential as a valuable resource in clinical radiotherapy. The research [14] delves into the realm of robustness analysis as a means to provide a more consistent framework applicable across various treatment techniques and modalities. This framework aims to address the uncertainties inherent in treatment planning and delivery, offering a standardized approach for evaluating and reporting plans. By identifying critical elements and dosimetric effects of uncertainties, robustness analysis seeks to enhance the reliability of plan evaluation, particularly in multi-institutional clinical trials. This approach holds the promise of promoting more accurate and consistent reporting of treatment outcomes, ultimately benefiting patients through more reliable radiotherapy. The research [15] presents an innovative concept of motion uncertainty, utilizing PDF to characterize motion caused by respiration. This concept is subsequently applied to construct a robust optimization framework for IMRT. Actual patient data is integrated into the analysis to assess the reliability of the generated solutions, using a clinical case of lung cancer as an illustrative example. The results are enlightening, showing that the robust solution effectively mitigates the under-dosing of the tumor compared to the nominal solution, particularly in worst-case scenarios. Furthermore, the robust approach showcases a significant decrease in the total dose administered to the primary organ at risk, specifically, the left lung. This observation underscores the capacity of this robust framework to enhance the optimization of radiotherapy by achieving an equilibrium between safeguarding healthy tissues and guaranteeing sufficient tumor dose delivery, a pivotal facet of radiotherapy planning.

In the paper [16], the emphasis lies on assessing the dosimetric effectiveness of robust optimization within the realm of helical IMRT for localized prostate cancer. The study involves a comparison of two distinct planning strategies: robust optimization and the conventional approach utilizing a planning target volume PTV margin. The evaluation considers various factors, including setup uncertainty and anatomical changes, both of which significantly impact treatment outcomes. The results suggest that robust plans exhibit potential benefits, including higher target coverage and lower organ-at-risk (OAR) doses, especially when perturbed scenarios are considered. However, the study also highlights the complexity of assessing robustness, particularly in the presence of anatomical changes. The article [17] introduces a ground-breaking concept of incorporating time-dependent uncertainty sets into robust optimization. This advancement tackles a prevalent issue in medical decision-making, particularly in situations where a patient's condition may evolve throughout the treatment process. In IMRT, such changes in cell oxygenation can directly impact the body's response to radiation treatment. The proposed framework offers a versatile approach to adapt to evolving uncertainties by modelling temporal changes within a cone structure, yielding current uncertainty sets at each treatment stage. The conic robust two-stage linear problems presented in this study cover a range of radiotherapy scenarios, and the clinical application of this approach is demonstrated in a prostate cancer case. The time-dependent robust approach is proven to improve tumor control over the course of treatment without introducing additional risks compared to established clinical methods. Furthermore, the research offers valuable insights into the timing of observations, maximizing the informational value for intermediate diagnostics. This innovative approach has implications not only in clinical settings but also in various applications, including maintenance scheduling.

### 3. Model Uncertainty

The objective of the motion PDF technique is to establish an accurate dose distribution by convolving it with an approximated PDF, thus addressing the problem of motion producing dose dispersion throughout radiotherapy [18]. However, this method requires prior knowledge of the expected motion pattern during treatment. If the actual motion pattern differs significantly from the assumed one, convolving with an optimized dose distribution for a different PDF can result in an uneven dose distribution with healthy and affected regions. As a result, a strategy to reduce treatment-related PDF uncertainty is required. Our conceptual framework is built on a finite set $X$ that represents the various components of the respiratory cycle. A PDF of motion is a nonnegative real function $f: X \to R$ that satisfies $\sum_{x \in X} f(x) = 1$. We begin with a nominal PDF designated as p, which was obtained from data gathered over the planning phase. We postulate that the nominal PDF $p$ may differ from the real PDF $\tilde{p}$ inside a subset $U$ of the domain $X$ and that this deviation is likely to occur after treatment.

This deviation follows an inequality condition.

$$p(x) - \underline{p}(x) \leq \tilde{p}(x) \leq p(x) + \overline{p}(x) \qquad \forall\, x \in U \quad [1]$$

In this equation, the values of $\underline{p}$ and $\overline{p}$ are used to define the range within which the assumed and actual PDF should not deviate. To ensure simplicity, we assume that $p - \underline{p} \geq 0$ and $p + \overline{p} \leq 1$ when considering the disparities during treatment.

$$P_U = \begin{cases} \tilde{p} \in R^{|X|} : \tilde{p}(x) \in \left[p(x) - \underline{p}(x), p(x) + \overline{p}(x)\right] \forall\, x \in U; \\ \tilde{p}(x) = p(x) \forall x \in X \backslash U; \sum_{x \in X} \tilde{p}(x) = 1 \end{cases} \quad [2]$$

It's worth emphasizing that incorporating the set $U$ can be considered somewhat superfluous since its impact can be efficiently managed by configuring both $\overline{p}(x) = \underline{p}(x)$ as equal to 0 for $x \in X \backslash U$, excluding those that pertain to $U$. The upper and lower bounds on $\tilde{p}$, which help define the range of uncertainty, will be referred to as "error bars".

The robustness of a treatment plan can be evaluated by checking that all of the conditions in our formulation are met regardless of which pdf from the set $P_U$ is actualized. During the optimization process, the set $P_U$ comprises all pdfs that must be protected against. Simple linear "smoothness" requirements can be incorporated into the concept of $P_U$ to alleviate concerns about conservative techniques allowing for implausible and extremely oscillating PDFs. $\tilde{p}(x) - \tilde{p}(y)| \leq \in$ if $|x - y| \leq \delta$, with suitable values for $\in$ and $\delta$. The primary challenge here is to ensure that $P_U$ encompasses a wide enough range of PDF variations to account for realistic patient-specific breathing patterns while preventing an excessive margin that would sacrifice critical patient information.

## 4. Optimization

The optimization techniques employed in this research are detailed in this section with the pseudo code.

**Cuckoo Search Optimization**

The cuckoo bird's brood parasitic behaviour served as inspiration for the CSO technique, which was initially presented in a journal [19]. Cuckoos use this behaviour to ensure that their eggs hatch by host birds [20]. For optimization, researchers looked at this natural process and developed the CSO method. The terminology within the CSO algorithm is metaphorically associated with familiar concepts in general optimization [21]. When dealing with single-objective function challenges, a "nest" or an "egg" represents an individual solution. An individual (nest) may contain many solutions (i.e., several eggs) in the arena of multi-objective function challenges. However, the primary focus of this research is on challenges with a single objective function. A nested set represents the entire population of possible solutions. The idea of a foreign bird departing a nest, representing the finding of a cuckoo's egg by a host bird, corresponds to the removal of an unsatisfactory solution. Conversely, the act of a cuckoo laying a new egg(s) in one or more nests represents the introduction of fresh solution(s) to the population. The CSO technique generates improved solutions using the following formula:

$$x_p^{t+1} = x_p^{(t)} + \otimes \; L\dot{e}vy(\lambda) \qquad [3]$$

Here, $x_p^{t+1}$ represents a fresh solution for a cuckoo labeled as $p$ acquired during a new iteration denoted as $t + 1$. This solution is derived from a prior solution, $x_p^{(t)}$, obtained in the preceding iteration $t$. To update these solutions, the Levy flight distribution algorithm called $L\acute{e}vy(\lambda)$ is employed. Here, $\lambda$ denotes the Lévy walk parameter, while $\alpha$ corresponds to the step size, which is determined by the scale of the particular issue being addressed. Furthermore, the symbol $\otimes$ represents element-wise multiplication. The Levy flight function approach allows for a stochastic walk with random step lengths derived from a Levy distribution. The Mantegna algorithm is typically used to estimate this distribution in the following manner:

$$L\dot{e}vy(\lambda) \sim \frac{u}{v^{-\lambda}} \qquad [4]$$

Where:

$$u \sim N(0, \sigma_u^2)$$
$$v \sim N(0, \sigma_v^2)$$
$$\sigma_u^2 = \frac{\Gamma(1+\lambda) * \sin\left(\frac{\pi\lambda}{2}\right)}{\Gamma\left(\frac{1+\lambda}{2}\right) * \lambda * 2^{\left(\frac{\lambda-1}{2}\right)}}$$
$$\sigma_v^2 = 1$$

With the Gamma function denoted as $\Gamma$. It's important to note that the parameter $\lambda$ falls within the range $1 < \lambda \leq 3$. The pseudo-code for CSO is given below.

PSEUDOCODE FOR CSO
*Define an objective function $f(X)$, where X represents the vector $X = (f(x_1, x_2, \ldots, x_d)^T$.*
*Initiate the population of host nests as $X_i (i = 1, 2, \ldots, n)$*

**While** $t < Max\_iterations$:
    Pick a cuckoo at random using Levy flights.
    Assess its quality or fitness denoted as $F_i$
    Randomly pick one nest among the $n$ nests (example $j$).
    **If** $F_i > F_j$,
        Replace $j$ with the new solution.
    **End if**
    Abandon a fraction ($pa$) of the less-fit nests and construct the new solution.
    Keep the best solutions.
    Determine the best solution by ranking them.
    Continue until the maximum number of iterations ($Max\_iterations$) is reached.
**End While**

## Flower Pollination Algorithm

The characteristics of the pollination process, pollinator behaviour, and flower constancy can be distilled into a set of rules in order to better understand them [22]:

1. Pollen-carrying insects engage in Lévy flights, allowing for the possibility of biotic and cross-pollination to occur on a global scale.
2. Conversely, self and abiotic pollination are examples of local pollination techniques.
3. The possibility of reproduction between two flowers is related to how similar they are to one another, and this is what we mean by flower constancy.
4. A switch probability, represented as $p$, ranging from 0 to 1 affects both local and global pollination probability. Local pollination assumes significance in overall pollination activities, influenced by factors like physical proximity and wind. The fraction $p$ signifies the contribution of local pollination to the entire pollination process.

Flowers can produce billions of pollen gametes, and some plants can have dozens of flowers on a single plant. For the purpose of simplification, we presume that each plant has a single flower and that this flower produces just one gamete of pollen. So, we can think of a solution $x_i$ as a gamete of pollen or a flower. This simplification could be further developed in the future to account for situations with multiple pollen gametes or numerous flowers in multi-objective optimization challenges. Depending on these idealized features, we may create a flower-based algorithm called the FPA. The two main phases of this method are called global and local pollination, respectively [23].

In the process of global pollination, insects and other long-distance travelers carry flower pollen from one location to another. Pollination and spread of the best ideas, represented by $g_*$, are therefore ensured. This process can be represented mathematically as follows, factoring in the first rule and flower constancy:

$$x_i^{t+1} = x_i^t + L(x_i^t - g_*) \quad [5]$$

In this equation, $x_i^t$ represents solution vector $x_i$ at $t$ iteration, and $L$ signifies pollination strength, serving as a step size. To mimic the variable step lengths observed in insects, a Lévy flight mechanism is employed efficiently. Then derive $L > 0$ from a Lévy distribution with the form:

$$L \sim \frac{\lambda \Gamma(\lambda) \sin(\frac{\pi \lambda}{2})}{\pi} \frac{1}{S^{1+\lambda}}, (S \gg S_0 > 0) \quad [6]$$

Here, $\Gamma(\lambda)$ represents the standard gamma function. Flower constancy and local pollination (Rule 2) could be depicted as follows:

$$x_i^{t+1} = x_i^t + \epsilon \left( x_j^t - x_k^t \right) \quad [7]$$

Flowers of identical species tend to remain consistent in appearance from one location to the next, and this depiction of pollen ($x_j^t$ and $x_k^t$) does the same. Mathematically, a local random walk involves uniform distribution sampling within [0, 1] if $x_j^t$ and $x_k^t$ are from identical species or populations. Most pollination of flowers occurs on both the local and global levels. Flower patches situated nearby or flowers within relatively close proximity are more susceptible to undergo local pollination compared to those positioned farther away. To this end, we suggest a proximity probability, denoted by $p$, to toggle between local and global pollination (Rule 4) [24]. The pseudo-code for FPA is given below.

PSEUDOCODE FOR FPA
*The objective is to minimize or maximize the function $f(x)$, where $x$ represents a d-dimensional vector $x =$*

$(x_1, x_2, \ldots, x_d)$. *Here are the steps of the algorithm:*
1. To begin, create a population of n pollen gametes or flowers, each of which will have a different random solution.
2. Find the optimal solution, $g_*$, among these possibilities.
3. Choose a switch probability $p$ between zero and one.
4. Proceed with iterations as long as the number of iterations $t$ remains below the maximum allowable generations denoted as "MaxGeneration." For each of the $n$ flowers in the population:
   - *If a randomly generated number falls below $p$, execute global pollination:*
   - *Create a d-dimensional step vector L following a Lévy distribution.*
   - *Update the global position using the equation:* $x_i^{t+1} = x_i^t + L(g_* - x_i^t)$.
   - *If the randomly generated number exceeds or equals "p," perform local pollination:*
   - *Generate a random value "$\epsilon$" from a uniform distribution between [0, 1].*
   - *Randomly select two solutions, j and k from the population.*
   - *Update the local position with the formula:* $x_i^{t+1} = x_i^t + \epsilon\,(x_j^t - x_k^t)$.
   - *Evaluate the newly obtained solutions.*
   - *If these fresh solutions demonstrate superiority over their predecessors, substitute them and locate the current best solution, $g_*$*
5. Continue this process until the achieved the maximum generation.

**Bat Search Optimization**

The bat algorithm is a bio-inspired technique rooted in echolocation, where bats employ sonar waves for navigation. It is a straightforward yet highly effective optimization method [25]. This approach draws inspiration from microbats' echolocation mechanisms, which these tiny creatures employ extensively to locate prey, identify obstacles in dark environments, and navigate through tight spaces, like stone cracks. The process of globally searching for a solution involves the position and velocity of virtual microbats undergoing random movements. Here, the position, referred to as $x_i$, represents the current value of the solution, while the velocity, $v_i$, indicates the transition from the current solution to potentially better solutions. At each iteration, the best current solution is indicated by $x_*$. The exploration of solutions involves adjusting parameters like frequency (wavelength) $f_i$, pulse emission rate $r$, and loudness Ai for each iteration. The effectiveness of this approach in locating global solutions depends on the precise management of frequency or wavelength to regulate the behavior of virtual microbats and achieve an optimal equilibrium between exploration and exploitation [26]. The mathematical equations governing the updates of location and velocity for each microbat in the group are outlined below:

$$f_i = f_{min} + (f_{max} - f_{min})\beta \qquad [8]$$

$$V_i^t = V_i^{t-1} + (X_i^{t-1} - X_*)f_i \qquad [9]$$

$$X_i^t = X_i^{t-1} + V_i^t \qquad [10]$$

Here, $\beta$ belongs to the range $[0,1]$ and it indicates the random vector from uniform distribution. The parameter $f_i$, signifying frequency (or wavelength), governs the rhythm and extent of the virtual bat's movement (both position and velocity) towards the local solution $x_*$ in each iteration and, ultimately, the best global solution once the objective is met. Additionally, the Bat Algorithm's efficiency is influenced by parameters like loudness and pulse emission rate. The mathematical expressions illustrating changes in sound value and pulse emission rate exhibit similarities, as illustrated below:

$$A_i^{t+1} = \alpha A_i^t \qquad [11]$$

$$r_i^{t+1} = r_i^0[1 - \exp(-\gamma t)] \qquad [12]$$

Where $0 < \alpha < 1 \; dan \; \gamma > 0$.

The Bat Algorithm operates under three key assumptions:
- All bats within the swarm employ echolocation for distance detection and the distinction between food and other objects.
- Bats navigate through random flight patterns, tuning their frequency (or wavelength) and pulse emission rate $(r)$ of sonar signals to determine subsequent positions and velocities. While the loudness value $A_i$ can fluctuate, it must remain within the range spanning from a high positive value $A_0$ to its minimum threshold, $A_{min}$.

- To gain a clearer understanding of the Bat Algorithm, the optimization approach is summarized in the pseudo-code below.

PSEUDOCODE FOR BSO
*Begin by initializing the bat population represented by $x_i$ and $v_i$ $(i = 1,2,....,n)$*
*Set the initial values for frequencies $f_i$, pulse rates $r_i$, and loudness $A_i$.*
**While** *(t < maxvalue), Proceed with the following steps:*
   *Generate novel solutions by adjusting the frequency using the formula:*
   $f_i = f_{min} + (f_{max} - f_{min})\beta$
   *Update the velocities and locations/solutions as follows:*
   $V_i^t = V_i^{t-1} + (X_i^{t-1} - X_*)f_i$
   $X_i^t = X_i^{t-1} + V_i^t$
   **If** *If $(rand > r_i)$, then*
     *Choose one solution from the good solutions available.*
     *Produce a local solution in the vicinity of the selected good solution $(x_*)$*
   **End if**
   *Produce a new solution through random flight.*
   **If** *$(rand < A_i \& f(x_i) < f(x_*)$, then*
     *Keep the new solutions.*
     *Enhance $r_i$ and decrease $A_i$.*
   **End If**
   *Rank the good solutions and identify the current good solutions*
**End While**

## 5. Result and Discussion

In this section, we analyze the results and discussion of IMRT on a tumor located in the lower left lung, considering the effects of respiratory uncertainty. Three optimization techniques, namely CSO, FPA, and BSO, were employed to optimize the treatment plan. The primary goal of the research is to guarantee that the correct dose was delivered to the tumour while minimizing the dose to healthy tissues, particularly the lung, and heart, which are prone to radiation-induced damage because of breathing fluctuation. To limit dose fluctuations, a constraint was set up on keeping a dose range between 72 Gy and 80 Gy.

In Figure 1, we see the Dose-Volume Histograms (DVH) for the tumor area after applying the three different optimization strategies. The red line indicates BSO, the blue line is FPA, and the black line is CSO outcome. Positive outcomes from all three optimization strategies suggest that the model effectively administered the prescribed radiation dose to the tumor while avoiding over- and under-dosing. Indicating that the three optimization strategies are adequate for guaranteeing the appropriate tumor dosage.

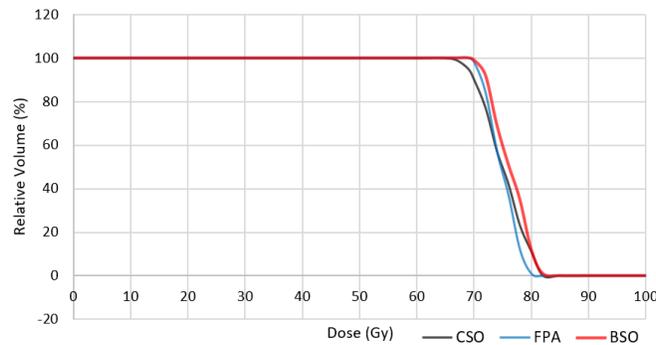

Fig. 1. DVH on tumor

Figure 2 depicts the DVH assessment for the lung area containing the tumor. The blue plot represents the dose distribution by BSO, the yellow plot represents FPA, and the red plot represents CSO. When the results of the three procedures are compared, it is clear that CSO delivered a minimal dosage to the lung. This finding suggests that CSO is the most successful approach for limiting radiation exposure to healthy lung tissue, hence reducing possible radiation-induced damage.

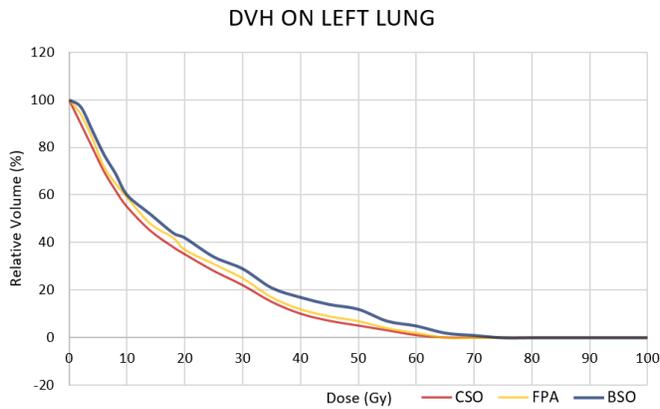

Fig. 2. DVH on Lung

Radiation can also affect the heart, hence Figure 3 displays the DVH analysis for that area as well. The dose delivered by the BSO is represented by the violet plot, FPA by the orange plot, and CSO by the blue plot. Similar to the lung region, CSO delivered the minimum radiation dose to the heart when compared to the other two optimization techniques. This result underscores the effectiveness of CSO in safeguarding the heart from excessive radiation exposure.

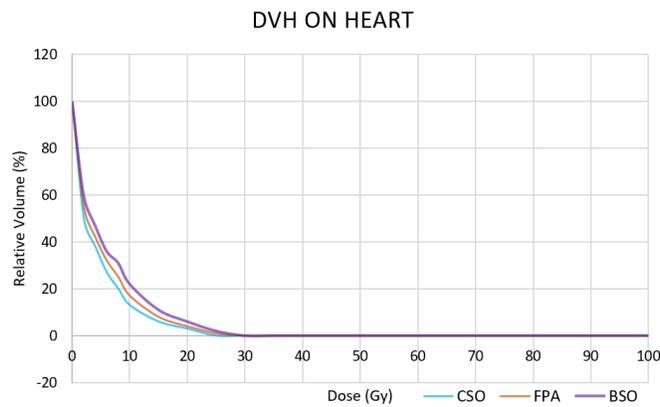

Fig. 3. DVH on Heart

Table 1 presents the doses delivered in Gy by the three optimization techniques to both the tumor and healthy organs. BSO and CSO delivered doses of 71 Gy and 70.32 Gy to the tumor, which are slightly lower than the minimum dose constraint of 72 Gy but within an acceptable range. However, FPA delivered a lower dose of 67.2 Gy to the tumor, indicating a deviation from the desired dose. CSO established its superiority in terms of healthy organs by delivering the least dosages of 29.4 Gy and 7.98 Gy to the lungs and heart. FPA supplied the second lowest dosage.

Table 1. Dose delivered to the tumor and other organs

| Optimization Method | CSO | FPA | BSO |
|---|---|---|---|
| Tumor | 70.32 | 67.2 | 71 |
| Lung | 29.4 | 30.554 | 35.924 |
| Heart | 7.98 | 10.25 | 12.57 |

In conclusion, the results show that CSO is the best optimization strategy for IMRT in the circumstance of respiratory uncertainty. It efficiently delivers a sufficient dose to the tumor while minimizing radiation exposure to critical organs, particularly the lungs, and heart.

## 6. Conclusion

In this study, we investigated the use of three bio-inspired optimization techniques—CSO, FPA, and BSO—to overcome the issues of optimizing radiotherapy under the conditions of respiratory uncertainty. Our main goal was to obtain a balanced dose distribution, assuring that the tumor received a sufficient amount of radiation while minimizing the dose given to healthy organs. Several significant inferences can be drawn from the results of this research and analysis. CSO and BSO were shown to be the most effective of the three bio-inspired approaches in terms of providing an adequate dosage to the tumor area. These methods were successful in confining the radiation dose to the tumor, which is essential for treating cancer. In particular, CSO proved to be the best method for radiation therapy planning when respiratory uncertainty was present. It not only provided the necessary dose to the tumor, but it also reduced radiation exposure to vital healthy organs like the lung and heart. This function is crucial in protecting surrounding tissues from damage caused by radiation. As we conclude this study, it is worth noting that while we focused on addressing respiration uncertainty, other sources of uncertainty in radiation therapy planning may exist. In future research, it is imperative to identify and tackle these additional uncertainties using bio-inspired or hybrid algorithms. By broadening the scope of optimization techniques and considering various uncertainties, we can further refine radiation therapy planning and enhance its precision and safety.